\documentclass[a4paper, 12pt, leqno]{amsart}
\usepackage{mathrsfs}
\usepackage{amsmath}
\usepackage{amscd}
\usepackage{amssymb}
\usepackage{amsthm}

\numberwithin{equation}{section}

\begin{document}
\title[Module Structure on Invariant Jacobians]
{Module Structure on Invariant Jacobians}

\author{ Nanhua XI}
\address{HUA Loo-Keng Key Laboratory of Mathematics and Institute of Mathematics\\
Chinese Academy of Sciences\\
Beijing, 100190\\
China } \email{nanhua@math.ac.cn}
\thanks{}

\maketitle

\def\Cal{\mathcal}
\def\bold{\mathbf}
\def\ca{\mathcal A}
\def\cdz{\mathcal D_0}
\def\cd{\mathcal D}
\def\cdo{\mathcal D_1}
\def\bold{\mathbf}
\def\l{\lambda}
\def\le{\leq}

Let $G$ be a connected semisimple algebraic group over an algebraically closed field $K$ of characteristic 0 and  $\rho$  a rational representation of $G$ on $V=K^n$. Then $G$ acts on the ring $A=K[x_1,x_2,...,x_n]$ of polynomial functions on $K^n$,  the representation is denoted by $\tau$.  The subspace $A_d$  of $A$ consisting of homogenous polynomials of degree $d$ is a $G$-submodule of $A$. For a polynomial $f$ in $A$, the Jacobian $J(f)$ is the subspace $A$ spanned by all the partial derivatives $\frac{\partial f}{\partial x_i}$.  When $G=SL_2(K)$, in 1985 Stephen Yau conjectured that if $J(f)$ is invariant, then the set of highest weights of $J(f)$ is a subset of the highest weights of $A_1$. (Convention: in this paper invariant means $G$-invariant.) Yau's conjecture arose from his study of singularities, see [NY, Y1-Y3, Y5].

 The conjecture was proved for the following cases: (a)  $n\le 5$, see [Y4], (b) $A_1$ being irreducible, see [SYY], (c) some special cases for $n=6,8$, see [Yu, YYZ]. For general $G$ and $f\in A_d$, Kempf showed that when $f=0$ is projectively smooth or $\rho$ is irreducible, if $J(f)$ is invariant and $d\ge 3$, then $f$ is invariant, see [K].  Kempf also showed in the same paper that if  $f$ is homogenous of degree greater than 2 and $J(f)$ is invariant, then there is an invariant homogeneous polynomial $g$ in $A$ with the same degree as of $f$ such that $J(f)=J(g)$. It is  clear that $J(f)$ is invariant if $f$ is an invariant, see [Y4, K]. The purpose of this paper is to give a positive answer to Yau's conjecture for arbitrary $G$. More precisely, we have the following result.

 \medskip

 \noindent{\bf Theorem.} Keep the notations above. Let $e_1,...,e_n$ be the standard basis of $V$ so that $x_i(e_j)=\delta_{ij}$. Then 

\noindent  (a) The linear map $A\otimes V\to A$ defined by $f\otimes e_i\to \frac{\partial f}{\partial x_i}$ is a homomorphism of $G$-modules.

\noindent (b) If $f\in A$ is invariant, then $J(f)$ is a quotient module of $V$.

\noindent (c) If $f\in A$ is homogenous of degree greater than 2 and $J(f)$ is invariant, then $J(f)$ is a quotient module of $V$.

\vskip3mm

{\sl Remark:} Since $A_1$ is the dual module of $V$, so that for $G=SL_2(K)$, $V$ is isomorphic $A_1$ as $G$-modules,  Yau's conjecture then follows from (c).

\vskip3mm

Proof. (a) Let $f=x_1^{a_1}\cdots x_n^{a_n}$. We need to show that the map sends $\tau(t)f\otimes \rho(e_i)$ to $\tau (t)(\frac{\partial f}{\partial x_i})$ for any $t$ in $G$. Let $\rho(t)(e_i)=\sum_{j=1}^na_{ji}(t)e_j$, then $\tau(t)(x_k)=\sum_{j=1}^na_{kj}(t^{-1})x_j$. So $\tau(t)f\otimes \rho(e_i)$ is sent to
$$\begin{array}{rl}&\sum_{j=1}^na_{ji}(t)\frac{\partial (\tau(t)f)}{\partial x_j}\\[3mm]
=&\sum_{j=1}^na_{ji}(t)\sum_{k=1}^n a_k a_{kj}(t^{-1})\tau(t)(x_1^{a_1}\cdots x_{k-1}^{a_{k-1}}x_k^{a_k-1}x_{k+1}^{a_{k+1}}\cdots x_n^{a_n})\\[3mm]
=&\sum_{k=1}^na_k(\sum_{j=1}^na_{kj}(t^{-1})a_{ji}(t))\tau(t)(x_1^{a_1}\cdots x_{k-1}^{a_{k-1}}x_k^{a_k-1}x_{k+1}^{a_{k+1}}\cdots x_n^{a_n})\\[3mm]
=&\sum_{k=1}^n\delta_{ki}a_k\tau(t)(x_1^{a_1}\cdots x_{k-1}^{a_{k-1}}x_k^{a_k-1}x_{k+1}^{a_{k+1}}\cdots x_n^{a_n})\\[3mm]
=&a_i\tau(t)(x_1^{a_1}\cdots x_{i-1}^{a_{i-1}}x_i^{a_i-1}x_{i+1}^{a_{i+1}}\cdots x_n^{a_n})\\[3mm]
=&\tau (t)(\frac{\partial f}{\partial x_i}).\end{array}$$ Hence the linear map is a hommorphism of $G$-modules.

(b) follows from (a).  Assume that $f\in A$ is homogenous of degree greater than 2 and $J(f)$ is invariant. According [K, Theorem 13], there is an invariant homogeneous polynomial $g$ in $A$ with the same degree as of $f$ such that $J(f)=J(g)$, thus (c) follows from (b). The theorem is proved.

\medskip

For the case of $SL_2(K)$, using the setup of [Y] (see also [SYY, Yu, YYZ]), one can show easily that the linear map $A\otimes A_1\to A$, $f\otimes x_i\to (-1)^i\frac{\partial f}{\partial x_{i'}}$ is a homomorphism of $SL_2(K)$-modules (or more precisely,  homomorphism of $sl_2(K)$-modules since the Lie algebra $sl_2(K)$ is used in [Y]), here $i'$ is determined by the following conditions (1) $x_{i'}$ is in the submodule generated by $x_i$ and (2) the weight of $x_{i'}$ is opposite to that of $x_i$. Combining [K, Theorem 13], this argument also proves Yau's conjecture.

\vskip5mm

\noindent{\bf Acknowledgement:} I thank Stephen Yau for talking me his interesting conjecture during a lunch after his colloquium talk at the Institute of Mathematics, Chinese Academy of Sciences.


\begin{thebibliography}{99}


\bibitem[K]{K} Gorge R. Kempf. {\sl Jacobians and invariants}. Invent Math, 112 (1993), 315¨C321.

\bibitem [NY]{NY} John N. Mather and Stephen S.-T. Yau. {\sl Classification of isolated hyper-
surface singularities by their moduli algebras}. Inventiones
Mathematicae 69 (1982), 243-251.

\bibitem[SYY]{SYY}  John Sampson,
Stephen S.-T. Yau, and Yung Yu. {\sl Classification of gradient space as
sl(2, C)-module I}. Amer. J. Math, 114 (1992), 1147-1161.

\bibitem[Y1]{Y1} Stephen S.-T. Yau. {\sl Continuous family of finite-dimensional
representations of a solvable Lie algebra arising from
singularities.} Proceedings of the National Academy of Sciences of
the United States of America, Vol. 80, No. 24 (1983), 7694-7696.

\bibitem[Y2]{Y2} Stephen S.-T. Yau. {\sl Solvable Lie algebras and generalized Cartan
matrices arising from isolated singularities}. Mathematische
Zeitschrift 191 (1986), 489-506.

\bibitem[Y3]{Y3} Stephen S.-T. Yau. {\sl Singularities defined by $sl(2, \mathbf{C})$ invariant
polynomials and solvability of Lie algebras arising from isolated
singularities}. American Journal of Mathematics, Vol. 108, No. 5 (1986),
1215-1239.

\bibitem[Y4]{Y4} Stephen S.-T. Yau. {\sl Classification of Jacobian
ideals invariant by sl(2, C) actions}. Mem. Amer. Math. Soc.,
72 (1988), pp.180.

\bibitem[Y5]{Y5} Stephen S.-T. Yau. {\sl Solvability of Lie algebras
arising from isolated singularities and nonisolatedness of
singularities defined by $sl(2, \mathbf{C})$ invariant polynomials.} American
Journal of Mathematics, Vol. 113, No.5 (1991),773-778.

\bibitem[YYZ]{YYZ} Stephen
S.-T. Yau, Yung Yu, and HuaiQing Zuo. {\sl Classification of gradient
space of dimension 8 by a reducible $sl(2, \mathbf{C})$ action}. Science in
China Series A: Mathematics, 52 (2009), 2792-2828.

\bibitem[Yu]{Yu} Yung Yu. {\sl On Jacobian ideals inviariant by a reducible $sl(2, \mathbf{C})$ action}.
Trans. Amer. Math. Soc., Vol. 348, No. 7 (1996), 2579-2791.


\end{thebibliography}
\end{document}